\documentclass[11pt]{amsart}
\usepackage{amsfonts}
\newtheorem{statement}{\bf Statement}
\newtheorem{definition}{\bf Definition}
\newtheorem{note}{\it Note}

\newtheorem{corollary}{Corollary}
\newtheorem{example}{Example}
\newtheorem{lemma}{Lemma}

\begin{document}
\title{Existence of the Ehresmann connection on a manifold foliated by the locally free action of a commutative Lie group}
\author{Pyotr N. Ivanshin}
\address{N.G. Chebotarev's RIMM, Kazan State University, Kazan, 420007, Russia}
\email{Pyotr.Ivanshin@ksu.ru} 
\begin{abstract}
In this paper the author determines necessary and sufficient
conditions for existence of the Eif hresmann connection on a
manifold foliated by locally free action of the
commutative Lie group. Also here we describe structure of
$C_0(M)|_L$ for a leaf $L \subset M$ in case such a connection exists. Finally we give some results
on structure of the spectrum of the family of Schr\"odinger
operators related to the leaves of the foliation.
\end{abstract}
\maketitle

\section*{Introduction}
In this paper the author considers some special structures on 
foliated manifolds \cite{Mo, Tam}. The main problem with foliated
manifolds is that they generally do not possess a good transversal
set (analogue to the base of the fibre bundle \cite{KoNo}).
Nevertheless one can find a set whose features are close to that
of the base of the fibre bundle. For example we can consider a
so-called Ehresmann connection on the foliated manifold \cite{BH}.
It is not clear though when we have such a structure on the
manifold. The aim of this paper is to find out the  necessary and sufficient conditions on
thefoliated manifold for the Ehresmann connection to exist. Note that here we consider only 
manifolds with the foliation
given by the locally free action of the commutative group. Note also
that once this Ehresmann connection exists one can reach some interesting
results (e.g. \cite{BH, BH1, BH2, ZhCh}). The  existence question
is then important and is considered in the first part of the
paper. Clearly the existence of the total connected
transversal on $(M, F)$ is necessary for the existence of the Ehresmann
connection. The open question is then if this condition is sufficient. 
We consider here mostly foliations of codimension $1$ though there are some notes on the general case.
It is shown here that if the foliated manifold possess a compact
connected transversal then we can construct an invariant Ehresmann
connection on it. So in this situation existence of a total
connected transversal is also a sufficient condition. In case the
total connected transversal is not compact we also get some more or less
satisfactory results on the possibility of existence of the
Ehresmann connection.

So let the foliation $F$ on the manifold $M$ be generated by a
locally free action of a commutative Lie group $H$. Once we have
an invariant Ehresmann connection on $(M, F)$ we can prove that
there exists an almost everywhere continuous bijection of our
manifold to the product $P \times S$ where $P$ is the transversal
manifold and $S \subset H$ is a fundamental set for a factor-group
$H'=H/H_1$ (\cite{SW}, Statement 1 and Corollary 1). Using this fact we arrive to the considerations
in the rest of the paper. Namely, first we consider density of the
intersection of the leaf $L$ with the invariant transversal $P$ and get a
structure of the $C_0(M)|_L$ in general situation. The question on
the structure of this algebra was inspired by the paper
\cite{Cadet}, where the author considered  the case of the compact
group $H$. In the third section  we consider a spectral problem
for the set of operators parameterized by points on $P$ with
potentials --- functions of $C_0(M)|_L$. 

\section{Existence of the Ehresmann connection.}
\subsection{Construction of the Ehresmann
connection invariant with respect to the group action.}
Let $M$ be a manifold with the foliation $F$ of codimension $1$ generated by
the locally free action of the commutative Lie group
$\mathbb{R}^n$. Let $\mathrm{codim}(F)=1$. Let us consider closed
connected transversal  $P \subset M$.

Condition $(*)$. Assume the existence of a leaf $L \in F$, which
has more than one intersection point with the manifold $P$. Consider a continuous bijection $p: \mathbb{R} \to P$. Suppose for any 
 $h_i \in H_x=\{h \in H| h x \in P\}$ the existence of a differentiable mapping $h_i: \mathbb{R} \to H$ such that   $h_i(t_x)=e \in H$, $h(t) p(t) \in P$, here $p(t_x)=x$. Note that locally ---in a sufficiently small neighbourhood of the point x--- this mapping always exists. Let $M$ possess a riemannian metric $g$. Assume next that there exists an element $a_{i j} \in H$ for any pair of $h_i, h_j \in H_x$ such that $g|_{L(t)}(h_i(t)p(t), h_j (t)p(t))=a_{ij}$, here $a_{ij}$ does not deped on $t$. This means that we have such a translation of the set $L \bigcap P$ along transversal $P$ that is an isometry on $L$. Suppose this translation is not a turn. This again can always be done locally. Note that then $h_i(t)=h_i(t_x)+b_i(t)$, $h_i \in H_x$.

\begin{definition}
Let us call transversal $P$ invariant one if  $h P\bigcap P \neq
\emptyset$ implies $h P=P$.
\end{definition}

\begin{statement}
Let $(M, F)$ be a manifold with foliation of codimension $1$ given
by locally free action of the commutative group $H$. Let there
exist a connected closed transversal $P$, which meets the
condition $(*)$.

Then we can deform the action of the group $H$ so that the orbits
of the action of $H$ (leaves of the foliation $F$) preserve but
the transversal $P$ becomes invariant with respect to the action
of $H$.
\end{statement}

$\bullet$ Fix a point $x \in P$.

1) Condition $(*)$ (existence of the set $H_x$) implies the existence of $h_1 = \inf_{\|h\|}
H_x = \{h \in H| h x \in P\}$, where $\|\cdot\|$ is a standard
Euclidean norm in $\mathbb{R}^n$.

2) Let us consider next a map $p: \mathbb{R} \to M$, $p(\mathbb{R})=P$,
$p(0)=x$ and set for any $y \in P$ $\gamma(s h_1) (y)= h_1 (b_1(t))^{-1}
y$. 
To clarify the construction one must consider two coordinate
charts adapted to the foliations in the neighbourhoods of points
$x$ and $h_1 x$. Since $h_1$ is a local diffeomorphism there exists
$\varepsilon>0$ such that for $(x-\varepsilon, x +\varepsilon)
\subset P$ we have $h_1 (x-\varepsilon, x +\varepsilon) \bigcap
P=\{h_1 x\}$. Also since codimension of $F$ equals   $1$, there
exists locally --- in a neighbourhood of $(x-\varepsilon, x
+\varepsilon)$ --- a continuous map   $h_t: (x-\varepsilon, x
+\varepsilon) \to H$. Condition $(*)$ (existence of $h_i(t)$) implies that the deformation can be
defined for all $P$. Note that $\forall s \in \mathbb{R} \setminus
\mathbb{Z}$ $s h_1 P \bigcap P =\emptyset$ since otherwise $h_1$
is not  an $\inf H_x$.

3) This deformed action is correctly defined. If there exists a
point $x \in M$ $h_1 x = x$ then by condition $(*)$ (absence of turns) for any $y \in
L_x $ $h_1 y= y$. This action is naturally continuous on $P \times
[0, \|h_1\|] \times \mathbb{R}^{n-1}$  with respect to the
coordinate $s h_1 \in H$. It is continuous in a neighbourhood of
$0$ and since it can be continuously spread up to any point by
condition  $(*)$ (existence of $h_i(t)$, $h_i \in H_x$) it is continuous for any other  $t \in
\mathbb{R}$.

Let us apply this construction for all $h \in H_x$.

Note that each consequent application of the preceding
construction the dimension of the nondeformed group action lessens
by $1$. So we must show that this algorithm can not be applied
more than $n$ times. Condition $(*)$ (no turns) implies that we get a
subgroup of the group of isometries $\mathbb{R}^n \to
\mathbb{R}^n$ which does not contain rotations so it consists only
of  translations. Since $P$ is a closed submanifold of $M$ the
generators set of the considered subgroup is not more than $n$
\cite{RS4}.

One can prove the last fact otherwise: This subgroup of isometries
consists of the maps $Ax +b$, $A \in O(n)$, $b \in \mathbb{R}^n$.
Let there exist  $(A-I)^{-1}$, then $\|(A^n + \ldots +I)b\|\leq \|
(A-I)^{-1}\|\|A^{n-1} -I\|\|b\|\leq 2 \| (A-I)^{-1}\|\|b\|$ for
any $n \in \mathbb{Z}$. This violates the condition $(*)$ (the set of translations consists of isometries). Thus
there exists at least one coordinate for each transformation
$\alpha_i(x)=A_i x + b$, with respect to which
$\alpha_i$ is operator of the type  $\left(%
\begin{array}{cc}
  1 & 0 \\
  0 & B_{i} \\
\end{array}%
\right)x + \left(%
\begin{array}{c}
  b \\
  b_1 \\
\end{array}%
\right)$ ($b\neq 0$). Then there can be no more than one
translation in this direction since otherwise in the neighbourhood
of the point $0=\lim\limits_{k\to \infty} b m_k - b' n_k$ either
condition $(*)$ --- there are turns if the set  $B_i b_2$ is closed, or otherwise $P$ is not closed. The construction of this group
gives us its commutativity. Let us prove now that it is generated
by no more than $n$ elements. Commutativity implies that for all
$i, j$ $(A_i-I)b_j=(A_j-I)b_i$. Let  $A_i=
\left(%
\begin{array}{cc}
  1 & 0 \\
  0 & B_{i} \\
\end{array}%
\right)$, then for any $j$ one has $A_{j, 1 l}-\delta_{1 l} b_{i
l}=0$, so there are no more than $n$ different $A_j$, consequently  there
are no more than $n$ $b_j$.  Commutativity implies that $A_i$ are
diagonal matrices which together with the fact that the action of
these isometries conserves orientation on any subset of
$\mathbb{R}^n$ (since the foliation is generated by the action of
the orientable group $\mathbb{R}^n$) gives us only the identity
matrix for each $A_i$. Thus our group is a commutative finitely
generated subgroup of  $Iso(\mathbb{R}^n)$ without cyclic subgroup
or the translation group. $\triangleright$

Note that the condition $(*)$ (no turns) of the previous statement is
necessary for the existence of the invariant Ehresmann connection
if the set $L \bigcap P$ is finite and $L \cong H$. Otherwise we
get the translation $a$ along transversal that is not a shift so
the action of the $h$ which translates a point $x$ of $L \bigcap
P$ into the $a(x)$ will produce infinitely many points $L \bigcap
P \subset L$ which contradicts the assumption.

\begin{example}{\em
Let us consider the foliation on $M=\mathbb{R}^2 \times
\mathbb{S}^1$ as the natural fibre bundle with the standard leaf
$\mathbb{R}^2$. Let us take the transversal $P$ which overlap
torus $\mathbb{S}^1 \times \mathbb{S}^1$ naturally
inserted into $M$, finite  number of times but more than once (Seifert foliation).
Since the condition $(*)$ of the previous statement is not
fulfilled  (the desired pair is the pair of points lying on the
intersection of the transversal  $P$ and  $\mathbb{R}^2$), $P$
does not define Ehresmann connection invariant with respect to the
action of the additive group $\mathbb{R}^2$.}
\end{example}

\begin{note}
Suppose that we do not want to deform the action of the group.
Here we can define a sequence of transversals which may converge
pointwise to the transversal which will define the Ehresmann
connection. Let us take as in the beginning of the proof of the
previous statement $x \in P$, $h_1 \in H_x$. Set $P_1= h_1 P$ then
consider $P_2= (h_1 + 1/2 h(t)) P_1$ where  $h(0) = h_{12}(0)$,
for $x \in P$ one has $h^{-1}_{12}(x) x \in h_1 P$ and $h(t_0)=e$
for $P(t_0)= h_1 x$, at the same time $[0, 1] h^{-1}_{12} \bigcap
P_2 = \emptyset$. And so on. The sequence of transversals
$(P_n)_{n \in \mathbb{N}}$ converges  pointwise to an invariant
transversal if and only if the infinite sum $\sum\limits_{i \in
\mathbb{N}} a_i$ where  $a_i= h_{i, i+1}(x)$ converges on any leaf
of the foliation $F$.
\end{note}

\begin{example}{\em
Consider the foliation on  $\mathbb{S}^1 \times \mathbb{R}$
generated by the images of the lines , $l_b$: $y=a x +b$, (the
constant $a\neq 0$), $b \in \mathbb{R}$ under the natural map
$\mathbb{R}^2 \to \mathbb{S}^1 \times \mathbb{R}$, $(x, y) \mapsto
([x], y)$. Let $P$ be the transversal which is defined in the
coordinate chart  $U \cong \mathbb{R} \times (-1/2, 1/2)$, $x \in
\mathbb{R}$, $y \in (-1/2, 1/2)$ as follows:
$$
y= \left\{
\begin{array}{cc}
 0 & x\leq 0; \\
 x &  x \in (0, 1/4];\\
 -x +1/2 & x \in (1/4, 1/2]; \\
 0 & x>1/2
\end{array} \right.
$$
If we take this transversal as the first member of the sequence
described i the previous not then there exists the limit
transversal of this sequence which is given by the constant
equality  $y=0$.}
\end{example}

\subsection{Properties of invariant transversals.}
Let us consider the connected transversal $P \subset M$ of the
foliation $F$ on $M$ generated by the locally free action of the
commutative group $H$. Let u denote the isotropy group of the set
$P$ by $H_P=\{h \in H| \forall x \in P, h x \in P\}$.

\begin{lemma}
There is no path that lies in $H_P$ and connects $x$ and $y$ for
any pair $x, y \in H_P$, $x \neq y$.
\end{lemma}

$\bullet$ Assume the contrary. Let there exists a path $\gamma:[0, 1]
\to H$, $\gamma[0, 1] \in H$. Then for arbitrary point $x \in P$,
$(\gamma[0, 1]) x \subset P$. Thus $P$ is not a transversal. This
contradiction implies the result. $\triangleright$

\begin{statement}
If $P$ is an invariant transversal then for any point $x \in M$ the isotropy group
$H_x \subset H_P$.
\end{statement}

$\bullet$ Again assume the contrary. Then in a neighbourhood
$U(x)$ of the point $x \in P$ which does not meet the condition of
the statement for any $h_x \in I_x$ $h_x U(x) \not\in P$, but at
the same time $h_x U(x) \bigcap P \neq \emptyset$. So $P$ does not
define a connection. $\triangleright$

Let $P$ be an invariant transversal on the foliated manifold $(M,
F)$. Consider the distribution of $TM$ tangent to the manifold
$P$. Let us transfer this distribution with the help of the action
of $H$ to the whole manifold $M$. The final distribution is smooth
one since $P$ is a smooth submanifold of $M$ and the action of $H$
is also smooth. Let us define the Ehresmann connection on $(M, F)$ as
follows: we set $\Pi(t, t')= s(t) w(t')$ for the vertical curve
$v: [0, 1] \to L$, here $v(t)=s(t) v(0)$, $s: [0, 1] \to H$,
$s(0)=e$, $s(1) v(0)=v(1)$, and the horizontal curve $w: [0, 1]
\to M$, $w= h w_0$, here $w_0 \in P$. Thus we obtain the rectangle
$\Pi$, unique due to the commutativity of the group $H$.

Consider then the set $\Gamma$ of the invariant transversals on
$M$ with isotropy groups $H_{\gamma}$, $\gamma \in \Gamma$. Each
transversal of $\Gamma$ gives rise to an Ehresmann
connection. Let $\gamma, \gamma' \in \Gamma$. We shall now define
the map $h: \mathbb{R} \to H$, $h(t) \gamma(t)=\gamma'(t)$.
Consider a leaf $L \subset M$ and a pair of points $x \in \gamma$,
$x' \in \gamma'$, $x, x' \in L$. There exists $h' \in H$ such that
$x h' =x'$. Set $h : \mathbb{R} \to \mathbb{R}^n$ as follows: Put
$h(0)=h'$. Note then that in some neighbourhood of the curve
$l[0, 1] x$ $l:[0, 1] \to H$, $l(0)=e$, $l(1)=h'$ one can define a
horisontal translation of this curve along transversals $\gamma$ ¨
$\gamma'$. It suffices now to show that this map can be extended to
the complete transversals. There are several possibilities:

1) The ends of the translated curve have limit points $x_1, x_2$
on the both transversals. Let  $x_1, x_2$ lie on the same leaf
$L$. Then there either exists  $\lim\limits_{t \to x_1}h(t)$ or
not. In the first case set   $h(x_1):=\lim\limits_{t \to x_1}h(t)$
and extend the translation by the rule described above. Let now
$h(t)$ be unbounded as $t \to x_1$. Since the transversal $\gamma$
is complete there exists $x^{'}_1 \in \gamma \bigcap L$. Let us
consider  $\mathrm{Sat}(U(x^{'}_1)$. There exists $t \in U'(x_1)$
such that $\gamma(t)= h(t) x''$ for some $x'' \in U(x^{'}_1)$.
Hence $x_1=h(t) x^{'}_1$ which contradicts the assumption since in
the neighbourhood of the singular leaf one can translate one
transversal into another with the help of the continuous bounded
map. Now let $x_1, x_2$ lie on the different leaves. Then by
completeness of the both transversals there exists an interval
intersecting both the leaf passing through the limit point which
lies on the other transversal as well as all the leaves near it.
So the transversals are not invariant.

2) There are no limit points on both transversals. The
construction is finished and the map $h$ is correctly defined.

3) The last possibility --- there exists a limit point $x$ only on
one transversal. Consider $\mathrm{Sat}((x-\varepsilon,
x+\varepsilon))$. The point $y \in \gamma' \bigcap L_x$ since the
second transversal is complete by assumption. Consider then
$\mathrm{Sat}(y-\varepsilon, y+\varepsilon)$. Since there are no
limit points on $\gamma'$ this curve is infinitely close to the
leaf $L_x$ thus intersecting any adjacent leaf on the one or the
other side of the limit leaf $L_x$. Now by considering a leaf from
$\mathrm{Sat}(y-\varepsilon, y+\varepsilon)$, we arrive to the
contradiction with the invariance of $\gamma'$.

Note that for a fixed $s \in [0, 1]$, $sh(t) \gamma(t) \in \Gamma$.
Fix a transversal $\gamma_0 \in \Gamma$, then for each pair
$\gamma, \gamma' \in \Gamma$ one can define $\gamma + \gamma' \in
\Gamma$. For a fixed point $x \in \gamma_0$ there exist uniquely
defined maps $h, h': \mathbb{R} \to H$, $\gamma(t)=h(t)
\gamma_0(t)$, $\gamma'(t)=h'(t) \gamma_0(t)$. Consider now
$(\gamma+\gamma')(t)=(h(t)+h'(t))\gamma_0(t)$. Isotropy group of
the last transversal naturally coincides with $H_{\gamma}
H_{\gamma'}$. Thus one can define addition of the two transversals
along the path connecting it.

The same operation can be defined for two Ehresmann connections
generated by two invariant transversals. So let $\nabla$,
$\nabla'$ be two Ehresmann connections invariant with respect to
the action of the group  $H$. Fix a point $x_0 \in M$. Now let us
consider transversals $P$ and $P'$ horizontal with respect to
$\nabla$ and $\nabla'$ respectively. Let $P$ and $P'$ pass through
the fixed point $x_0$. Let us construct $P''=P+P'$ along zero path
$\gamma([0, 1])=x_0$. Thus there is a structure of the additive
group on the set of invariant Ehresmann connections.

\begin{example}{\em
Let us consider different flows on the torus  $\mathbb{T}^2$ as
the set of transversals. Let the leaves of the foliation be
meridians. Consider the set of transversals generated by flows $\Gamma
\supset \Gamma' \cong \mathbb{R}$, here the map $\Gamma' \to
\mathbb{R}$ is given by formula $\gamma \mapsto a$, $(a, 1)$ being
a vector of the flow $\gamma$.}
\end{example}

\begin{example}{\em
Consider the foliation on $M=\mathbb{R}^2\setminus\{0\}$,
generated by circles with common center $0$. Let us define an
action of the group $H=\mathbb{R}$ on $M$ as follows: $\mathbb{R}$
$R: \mathbb{R} \times \mathbb{R}^2 \to \mathbb{R}^2$: $R_h x= e^{i
h/|x|} x$. Then there is no connection invariant with respect to
the action of  $H$, but one can deform given action of
$\mathbb{R}$ as follows: $R^{'}_h x= e^{i h} x$. Let us show by
reductio ad absurdum that there is no invariant Ehresmann
connection on $M$ with the foliation generated by undeformed
action of $H$. Assume the existence of such connection, then the
tangent vector at a point  $x \in \mathbb{R}^2$ to the horizontal
transversal (which exists since the considered foliation is the
foliation of codimension equal to  $1$) under the action of the
element from isotropy group  $H_x\cong \mathbb{Z}$ of the leaf
$L_x \ni x$ transforms to the collinear one. This is impossible
since $\frac{d}{d r} e^{i h/r}=-\frac{1}{r^2} e^{i h/r} \neq 0$.}
\end{example}

Now let $\mathrm{codim} F >1$ and consider global connected
transversal $N \subset M$. Let us then consider a set of curves
$\Phi$ of $N$. Note that $\forall \phi \in \Phi$ the existence of
$h' \phi \bigcap \phi \neq \emptyset$ implies as in the proof of
the first statement the existence of the map $h: \mathbb{R} \to
H$, $h(t) \phi \subset N$. Let $h(t) \phi \in \Phi$. Let also the
condition $(*)$ be true for any curve of $\Phi$. Thus
either $h' h(t) \phi=\phi$ or $(h' \phi) \bigcap \phi$ is
discreet.  Consider the  subset $H^c_{\phi_0} \subset H_{\phi_0}$
defined for a curve $\phi_0 \in \Phi$. We deform the action of $H$
as in the proof of the first statement to get a connection on
$\mathrm{Sat}(\phi_0)$, invariant under the action of $h \in
H^c_{\phi_0}$. Let us then consider the set of curves $\Phi_0=
\bigcup\limits_{n \in \mathbb{Z}} h_1^n \phi$, where $h_1 =
\inf\lim_{\|\|} \{H^d_{\phi_0}\}$ and deform the action of $H$ on
$\mathrm{Sat}(\Phi_0)$.

1) Assume first that $h_1 d\gamma(x) \in T_{h_1 x} N$ for
$x=\gamma(0) \in N$ (this can always be done). Note that if there
exists a neighbourhood of identity $U(e) \in H$ such that
$U(e)\setminus\{e\} N \bigcap N = \emptyset$ then $\dim M
\setminus \overline{\Phi_0} = \dim M$.

2) Now we find ourselves in the conditions of the third statement.

3) If $\overline{\mathrm{Sat(\Phi_0)}}=M$ then the algorithm of
construction is finished, otherwise we must consider another curve
$\phi' \in \Phi$ and continue the construction for it.

1) implies that the dimension of the final distribution of $TM$ equals $p =\dim N$.

The construction algorithm implies that one can construct the
invariant connection on $M$ if there exists at least one
connection on $(M, F)$ in which any two leaves can be connected by
horisontal line. Note that the invariant connection may not be
integrable.

\subsection{Analysis of conditions of Statement 1.}

Now let transversal $P$ be not a closed subset of $M$. Then if
$\gamma \bigcap L$ is dense in $L' \subset L$, $L'=\{x_0\} H_1$,
for $H_1$ being a subgroup of $H$ it is sometimes possible to get
an invariant connection on $M$. Consider for $\gamma(0)=x_0$,
$t_1=\inf\{t \in [0, \infty)| \gamma(t) \in L(0)\}$, then apply
the algorithm to $h=\gamma(t_1)\gamma(0)^{-1}$. Note that in this
case  $\gamma$ generates the distribution on $TM$. We can apply
Statement 1 to $H \cong \mathbb{R}$, $P \cong L$. Note that the
translation $f : L \to L$ of the leaf $L$ along $\gamma$ defined
with the help of the first intersection point of $\gamma \bigcap
L$ with respect to a fixed point $x_0 \in \gamma \bigcap L$ can be
made an isometry since $\forall n \in \mathbb{Z}$, $f^n(x_0) \neq
x_0$. The same can be done in case the  leaf  intersecting $\gamma$
is dense or $\overline{\gamma\bigcap L}_{L}=L$ and
$\overline{\gamma\bigcap L}_{\gamma}=\gamma$. This means that we pass to
the isometry group on $L$. Note that the results above hold true
only under conditions $(*)$ for each of considered maps.

Now let $\gamma$ be not a closed subset of $M$ and not everywhere
dense. Consider $x_0 \in \overline{\gamma} \setminus \gamma$. Let
us try to deform action of the group according to the algorithm of
the first statement for $x_n \in \gamma \bigcap
B_{\frac{1}{n}}(x_0)$. Note that $h_1(n) \to 0$ as $n \to \infty$.
Then a point $x \in \omega$ after the deformation of the
group action becomes unaccessible from any other point of the
 leaf  $L$. Thus $\gamma$ can not be an invariant transversal.

\begin{note}
Let the transversal $\gamma$ be either not closed submanifold of
$M$ or violate condition $(*)$ (not an isometry). Then the method described in the
paper \cite{Nov} gives us a compact connected transversal $S$ on
the foliated manifold $M$. Note that later on we shall try to
avoid this method in order  either to deform existing non-compact
transversal or to get the compact one with some special
properties.
\end{note}

If $\omega(\gamma(0))$ being the $\omega$-limit set of $\gamma$ is
closed one-dimensional submanifold of $M$ then we can consider it
as new $\gamma$.

\begin{statement}
Let there exist a riemannian metric on $M$. Suppose for $u \in F \subset TM$, $\|u\|=1$, $(v, u)\leq \alpha <1$
for any $v \in TM$, $\|v\|=1$, $v$ being a tangent vector to
$\gamma$. Then if for a fixed leaf $L$ there exists an isolated limit point $x \in
\overline{\gamma} \setminus{\gamma}|_{L}$ then

{\em 1)} There exists a limit set $\omega$ of the transversal $\gamma$
which is homeomorphic to $\mathbb{S}^1$.

{\em 2)} $\omega$ is a total transversal of the foliation
$F$ if there exists at least one compact  leaf  intersecting  $\omega$.
\end{statement}

$\bullet$ 1) Let us first introduce a parametrization on
$\gamma(t)$ such that $x_k=\gamma(t + k)$. The set of functions
$f^k$ with graphs $\gamma|_{\gamma \bigcap x_k \neq \emptyset}$ is
uniform continuous in the charts adapted to the foliation in a
neighbourhood of $x$. Hence there locally exists a limit curve
$\omega(x)$ transversal to $F$. The curve $\omega$ is closed since
otherwise $x$ is not an isolated limit point. Sliding along
$\gamma(t)$, $t \in [0, 1]$ by compactness of this interval we get
a limit set $\omega$ (maybe not unique) which is naturally
homeomorphic to $\mathbb{S}^1$. The set $\omega$ is transversal to
the leaves of the foliation $F$ by construction.

2) Compactness of $\omega$ implies the existence of leaves which
lie in a neighbourhood of $\omega$. Let us show that there are no
other leaves. First consider $T=\inf\{ t \in \mathbb{R}| \forall
s> t \in L, L_s \bigcap \omega \neq \emptyset\}$. Since
$\mathrm{Sat}(\omega)$ is open $\gamma(T) \in L$, $L \bigcap \omega =
\emptyset$. Then $(T, +\infty)$ is covered by the infinite set of
intervals of the type $[a_n, a_{n+1})$, $n \in \mathbb{Z}$, for which
we have  $a_n \to T$ as $n \to - \infty$ and $a_n \to + \infty$ as $n
\to +\infty$. Thus the topological structure of the arbitrary
neighbourhood of the limit point $\gamma(T)$ is non-trivial which
contradicts with the definition of the foliation.$\triangleright$

\begin{note}
One can show compactness of the $\omega$-limit set of the
transversal $P$ in case $x$ is not an isolated limit point using
the same methods as in the proof of the closeness of $\omega(x)$.
At the same time little can be said on the dimension of this limit
set.
\end{note}

\begin{example}{\em
Assume that the foliation on torus
$\mathbb{T}^2=\mathbb{R}^2/\mathbb{Z}^2$ is given by parallels.
Consider the integral curve of an irrational flow as the
transversal $P$. Then for any point $x \in \mathbb{T}^2$
$\omega(x)=\mathbb{T}^2$.}
\end{example}

Let us now clarify the condition of the second part of the
previous statement by the following examples:

\begin{example}{\em
Let a foliation on $(-1, 1) \times \mathbb{S}^1$ be given by
gluing two Reeb foliations on $(0, 1) \times \mathbb{S}^1$. Then
there exists a complete transversal diffeomorphic to $\mathbb{R}$
which converges to $\{-1/2\} \times \mathbb{S}^1$ and $\{1/2\}
\times \mathbb{S}^1$. But neither of the limit sets is total
transversal since they do not intersect the leaf $\{0\}\times
\mathbb{S}^1$.}
\end{example}

\begin{example}{\em
Consider Klein bottle as foliated manifold $M$. Define the
foliation on it by meridians between two of which we put the Reeb
foliation on $(0, 1) \times \mathbb{S}^1$. Again any limit set of
any complete transversal does not intersect two boundary leaves.}
\end{example}

\begin{note}
Note that in case all the leaves passing through $\omega$ are
homeomorphic to $\mathbb{R}$ each leaf intersects $\omega$ only at
one point. Assume that each leaf of the foliation is noncompact.
Then there exists a subgroup $H_1$ of $\mathbb{R}^n$, homeomorphic
to $\mathbb{R}$ such that for any  $h \in H_1$ $h \omega \bigcap
\omega=\emptyset$.
\end{note}

$\bullet$ The following statements are true due to the facts proved above
together with the construction from the next statement.

1) There exists a connection on $\mathrm{Sat}(\omega)$ invariant
with respect to the action of the group $H$.
$\mathrm{Sat}(\omega)$ is an open set since each point $x$ of it
possess an open neighbourhood of the type $U(h) (\alpha_1,
\alpha_2)$, where $x= h \alpha$ for $\alpha \in (\alpha_1,
\alpha_2) \subset \omega$, $h \in U(h) \subset H$ and $U(h)$ being an
open set.

2) There exists a sequence $(h^n y) \to
L_x$, $n \to \infty$ for any $y \in \omega$, here $L_x$ is a leaf passing through the
point $x \in P$, which is a limit point of the previous statement
and $h^n \in  \{h \in H| h y \in \omega\}$. Thus $\forall h \in H$
$\exists n \in \mathbb{N}$, $\forall m> n$ $h^{m} p \bigcap \omega
=\emptyset$.

Hence for all $h \in
H\setminus\{e\}$
$h
\omega
\bigcap \omega = \emptyset$.
$\triangleright$

Note that the set $\mathrm{Sat}(P_1) \bigcap P_2$ is always an open subset of $P_2$ for
any two transversals $P_1$ and $P_2$ on the manifold with
foliation of codimension $1$.

\begin{statement}
Assume that condition $(*)$ (no turns part of it) is violated but the set $P$ is
closed hence homeomorphic to the circle $\mathbb{S}^1$. Then the
following statements hold true:

{\em 1)} The number of intersection points of the leaf $L'$ with $P$
between any fixed pair of the intersection points of the fixed
leaf $L$ with $P$ with respect to the translation of the pair of
points along transversal is constant.

{\em 2)} There exists a number $N \in \mathbb{N}$ such that for any $L
\in P$ the cardinality of the set $L \bigcap P$ is not more than
$N$ if the translation along transversal does not contain shifts.

{\em 3)} Again if the translation along transversal does not
contain shifts then there is no leaf intersecting $P$ infinitely
many times.
\end{statement}

$\bullet$ 1) It suffices to consider the construction of the first
statement with the only possible difficulty, i.e. there may not
exist $h \in H$, $h P=P$, since otherwise there again may exist $n
\in \mathbb{N}$, $h^n=e$. Since the translation along transversal
$P$ is a continuous map the proposition holds true.

2) The map generated by the map $x_i \to x_{i+1}$ generated by
translation along transversal can be lifted to the isometry  $U: H
\simeq \mathbb{R}^n \longleftrightarrow \mathbb{R}^n$; besides
since $U H=H$ and $U^n=Id$ we have $U \in O_n(\mathbb{R}^n)$ (the
set of orthogonal transformations of $\mathbb{R}^n$). Let us
consider the leaf $L$ intersecting $P$ minimal possible number of times (let
us denote this number by $n$) then each other leaf by 1) must
intersect $P$ by the set of cardinality multiple of $n$. Since the
transformation $U$ preserves while sliding along leaves the cyclic
group $<U(L')>$ is a subgroup of $<U(L)>$. Thus on any leaf we get
the action of the subgroup of the orthogonal translations and if
we do not assume the existence of the leaf that intersects $P$
maximal finite number of times then $P$ is not closed. This
contradiction completes the proof.

3) Let there exist a leaf $L$ such that $L \bigcap P$ is infinite.
Let us consider the set of orthogonal transformations $U(L, P)$ on
$\mathbb{R}^n$. Let this set be generated by the set of points $L \bigcap P$.
Since $U(n)$ is a compact set there exists a limit point for $U(L,
P)$, so again $P$ is not closed. $\triangleright$

\begin{note}
If $P$ is compact then the condition $(*)$ holds true.
\end{note}

$\bullet$ The proposition is obvious if all the leaves are
compact. Otherwise let us prove this fact for foliations of the
dimension equal to $1$ (the proof is similar in general case). Let
condition $(*)$ (no turns part) be satisfied. The same holds true in case there
exists a curvature since the operators $U(L)$ are unitary. Let the
leaf $L$ intersect $P$ infinite number of times then the only
obstruction while translating $[x_{i}, x_j] \subset L$ along $P$
is as follows: $x_{n+i} \to C$, $x_{j+n} \to \infty$ as $n \to
\infty$. But this contradicts the definition of the foliation
since then the set
 $(x_{n+i})_{n \in
\mathbb{N}}$ has a limit point on $P$ and this limit point is
infinitely close to its image under the translation by elements
from the neighbourhood of $C$. Or we can point out the existence
of $N \in \mathbb{Z}$ such that $x_{N+i}=x_j$. Note that there exists a correctly defined set $h_1(t)$, $t \in \mathbb{S}^1$ since the transversal is compact and this set is defined for a neighbourhood of any point of it. $\triangleright$

\begin{example}{\em
Let us describe the set $L \bigcap P \subset L$ and the set of its
transformations in some partial cases.

1. If $H=\mathbb{R}^1$, then the transformations preserving
discreet set of points and an orientation on $\mathbb{R}$ can be
only shifts.

2. If $H=\mathbb{R}^2$ then the following cases of the set $L
\bigcap P$ are possible:

a) $L \bigcap P$ is finite then the deformations which conserve it
are the turns on the angle $2\pi/k$, $k \in \mathbb{N}$.

b) $L \bigcap P$ is countable then the deformations can be only
turns by the angles:  $\pi/3$ ($\mathbb{R}^2$ by isosceles
hexagons), $2\pi/3$ (triangles), $\pi/2$ (squares), $\pi$
(stripes), along with the shifts in subsequent directions. The
invariant transversal then passes trough the centers of the given
figures. Note that by construction of the translation along $P$
the centers can not be points of $L \bigcap P$.}
\end{example}

\begin{corollary}
Let $M$ be a manifold with the foliation $F$ of codimension $1$,
generated by a locally free action of the commutative Lie group
$H$. Assume that there exists a compact transversal on $(M, F)$.
Then there exists an invariant Ehresmann connection on $(M, F)$.
\end{corollary}

$\bullet$ It suffices to show the existence of the one-dimensional
full transversal which meets the conditions of the previous
statement. First note that the map $U(L_M): \mathbb{R}^n \to
\mathbb{R}^n$, for $L_M$ being a leaf with the minimal number of
intersection points with $P$, naturally has a fixed point which
slides along $P$ to the neighbouring leaf $L$. Since $U(L_M)$ can
be uniquely extended to the neighbourhood of the leaf $L_M$ we get
the result. $\triangleright$

Condition $(*)$ (isometry part) can be analyzed in the following way: if $[x(t),
x'(t)] \to \infty$, then the unique non-trivial situation which
can not be improved with the help of the last method described
below is $x(t) \to x_0$, $x'(t) \to x^{'}_{0}$, where $x_0,
x^{'}_{0}$ belong to different leaves. Then there naturally exists
an element $h \in H$ such that the distance between $sh x_0$ and
$sh x^{'}_{0}$ infinitely increases as $s \in \mathbb{R}$. Let us
consider now the two-dimensional foliation on the set $sh \gamma$,
$s \in \mathbb{R}$. The only foliation such that the transversal
$\gamma$ can not be improved with the help of the method described
below is a Reeb foliation. Thus the transversal $P$ intersecting
leaves between $L$ and $L'$ is not compact or more precisely is
the union $(-\infty, a] \bigcup [b, +\infty)$. Then there also is
no closed transversal connecting $x_0$ with $x^{'}_{0}$ (otherwise
on this transversal exists a point which is invariant under the
action of $sh$ for all $s \in \mathbb{R}$). Hence the manifold $M$ splits
into two saturated sets as follows: $\mathrm{Sat}([a, b]) \bigcup
\mathrm{Sat}((-\infty, a] \bigcup [b, +\infty))$. Note that in the
latter case we can consider one of the ends $(-\infty, a]$ or $(b,
+\infty)$ as new transversal. Note also that there can be no more
than $|\pi_1(M)|$ of such obstacles. This is true since
there exists a loop which is a concatenation of the set $[a
+\varepsilon, b-\varepsilon]$ and an interval connecting point
$x_1=x(t_1)=a +\varepsilon$ with $x_1^{'}=x'(t_2)=b-\varepsilon$.
If this loop is contractible then there exists $h:[a +\varepsilon,
b-\varepsilon] \to H$ $h[a +\varepsilon, b-\varepsilon] [a
+\varepsilon, b-\varepsilon] ={pt}$, which is impossible.

Now let us consider two constructions that help us to deform a
transversal $\gamma$ so that it satisfy condition $(*)$ (isometry part).
Assume that $\gamma(t) \to \omega \subset L$, $t \to \infty$ and
that there exists  $t' \in \mathbb{R}$ such that $\forall t> t'$
$\omega(t) \not\in \gamma$. Note that in this case $\gamma((t',
+\infty))$ lies in the saturated neighbourhood of the leaf $L$.
Let us construct $\gamma'$ as follows: since there exists $t \in
\mathbb{R}$ such that $\gamma(t) \in \omega$, the set $H \gamma(t
-\varepsilon, t+\varepsilon) \bigcap \omega \neq \emptyset$ and
there exist $\delta>0$, $h \in H$ and $t \in \mathbb{R}$ such that
$\gamma(t) \in U_{\delta}(\omega)$. Note that for some $t_0>t'$
$\gamma([t_0, +\infty) \in \mathrm{Sat}(\gamma(t -\varepsilon,
t+\varepsilon))$. Let us construct $\gamma'$ by gluing of the set
$\gamma(-\infty, t)$ to the part of $\gamma$ converging to $L$
from the proper side, i.e. with $t'' \in \mathbb{R}$, $\gamma(t'')
\in \mathrm{Sat}(\gamma(t -\varepsilon, t+\varepsilon))$, $t''>
t'$. Thus we get a curve connecting $L$ with itself. The only
difficulty here is the possibility of self-intersection of the
deformed transversal $P'$. But by Note 2 there
then exists a transversal homeomorphic to $\mathbb{S}^1$.

The second case. Let $x \in M$ be a limit point of the transversal
$\gamma$. Also let there exist a coordinate chart such that for
$t_n \in \mathbb{R}$, $t_n \to \infty$, $|d\gamma(t_n)/d t| \to
\infty$, $(n \to \infty)$. Consider a number $C \in \mathbb{R}^+$.
Let us deform $\gamma$ as follows: consider a function $\phi_n:
\mathbb{R} \to H$, $\phi|_{(-\infty, t_n - \varepsilon)} = e$, $|d
\phi(t_n) \gamma (t_n)|< C$, $\phi|_{(t_n, +\infty)}=\phi(t_n)$ in
a neighbourhood of $t_n$. The  curve $\gamma'=\gamma
\prod\limits_{n \in \mathbb{N}} \phi_n$ then does not have the
same property as the curve $\gamma$ at the point $x$. Again there
may appear new points on $\gamma'$ with this property.

\section{Structure of $C_0(M)|_L$.}
Let $(M, F)$ be a foliated manifold with the foliation generated
by the locally  free action of a commutative group  $H$. Assume
that $H$ is a one-dimensional group. Assume also that $(M, F)$
possess an invariant integrable Ehresmann connection. Let $P$ be
an invariant transversal. Let us consider a continuous map $F: P
\to P$ generated by the action of $a \in H_P$. Since $C_0(M) \in
C_0 \times '([0, a])$ it suffices to describe structure of
$C_0(P)$. Let us consider a sequence $(d(f^n(x), f^n(y)))_{n \in
\mathbb{N}}$ for each pair of points $x, y \in P$. In case $0$ is
an accumulation point of this sequence then there are two
possibilities: 1) there exists a point $z \in P$ and subsequence
$(d(f^{n_i}(x), z)) \to 0$ as $n \to \infty$; 2) there is no
accumulation point. The second case is not interesting since then
it is possible to deform metric on $P$ so that $d(x^{n_i},
x^{n_j})\geq \varepsilon$ for some $\varepsilon>0$. In the first
case then the restriction of any function $g \in C_0(M)$ on the
leaf $L_x$ possess a subsequence $|g(a^{n_i} x) - g(z)| \to 0$.
Define a function $g_z$ on $\mathbb{R}^{+}= \bigcup\limits_{i \in
\mathbb{N}} [a^{n_i}, a^{n_i+1}]$ s a restriction of $g$ on the
set of intervals correspondent to the set of points $a^{n_i} x$.
The function $g_z$ behaves as in the third case of \cite{Lob}
where the component $f_1$ is bounded but does not belong to
$C_0(\mathbb{R})$.

Now let us describe the structure of $C_0(M)|_L$ for arbitrary
leaf $L \in F$. Fix $\varepsilon>0$. We shall consider a finite
number of limit points $z_k$, $k=\overline{1, . . n(\varepsilon,
K)}$ for any leaf $L$ on the compact subset of $P$. Then let us
put into consideration the set of functions $g_{z_i}$
approximating the given function $g \in C_0(M)|_{L}$ with
precision $\varepsilon$ on the compact subset $K$ of $P$ with
respect to the supremum norm. The function $g_{z_i}$ is defined on
the set of intervals $[a^{n_k}, a^{n_k+1}]$, here $a^{n_k} x \in
U_{\delta}(z_i)$, $\forall t \in [0, 1]$ $g_{z_i}(a^{n_k} t a x)=
g(a^{n_k} t a z_i)$. Note that $g_{z_i}$ is a bounded function
with discontinuities in the set of points $\bigcup\limits_{k \in
\mathbb{N}} \{a^{n_k} x\} \bigcup \{a^{n_k+1} x\}$. Now as in the
third case of \cite{Lob} we deform $g_{z_i}$ to $g^{'}_{z_i} \in
C_b(L_x)$, so that $\| (g-g^{'}_{z_i})|_{\bigcup\limits_{k \in
\mathbb{N}}[a^{n_k}, a^{n_k+1}]}\| < \varepsilon$. Consider now
$g_1= g - \sum\limits_{i=1}^{n(\varepsilon)}g^{'}_{z_i}$ as new
function and approximate it with precision $\varepsilon/2$ on a
larger compact $K'$ of $P$ and so on. As the result we get the
function $g' \in C_0(\mathbb{R}^+)$. Thus we get a sequence of the
finite sets of functions the behavior of which was described in
\cite{Lob} and which approximate the target function $g$ with any
precision on any compact subset of $P$. Hence
 we get an inclusion $C_0(M)|_L \subset C_0(L)+ \times_{\mathbb{Z}}
C([0, 1])$. The last space is a space with the norm defined by the
algorithm of approximation. That is this norm is as follows:
$\|f\|=\|f_0\|_0+\sum_{l \in
\mathbb{N}}\max\limits_{k_{l1}<i<k_{l2}}\{\|f_i\|_0\}$.

The only question that needs clarification is the question of independence of the
approximating point set $(z_i)\subset P$ on the target function.
To avoid this difficulty one must consider two functions $f_1,
f_2: M \to \mathbb{R}$ such that $\|f^{'}_1\|^2+\|f^{'}_2\|^2 \neq
0$. Let us then unite two sets of approximating points for these
functions into one. Let us show now that it turns possible to find
a function $g_i$ from the space constructed in the previous
paragraph such that $\|g -g_i\|< \varepsilon$ for any function $g
\in C_0(M)$ and  $\varepsilon>0$. Consider now $g_1 \in
C^{1}_0(M)$ such that $\|g -g_1\|<\varepsilon$, $g_1|_{M\setminus
K}=0$ on a compact $K \subset M$. Then there exists $C \in
\mathbb{R}^+$, $\|g^{'}_1\| \leq C$. Hence for all $x, y \in K$
$|g(x) -g(y)| \leq N |f(x)- f(y)|$, for $N=C/\max\{\|f^{'}_1\|,
\|f^{'}_2\|\} < \infty$. Thus there exists a function $g_i$ of a
desired space such that either $\|g_1-g_i\|_0 < \varepsilon$ or
$\|g -g_1\|<2 \varepsilon$.

\begin{statement}
{\em 1)} In the given decomposition of the algebra $C_0(M)|_L$ the first
component (that is $C_0(L)$) vanishes if and only if $L \bigcap P$
does not contain isolated points.

{\em 2)} The second part $\times_{\mathbb{Z}} C([0, 1])$ of the
mentioned decomposition is finite if and only if the set of limit
points of $L \bigcap P \subset P$ is finite.
\end{statement}

$\bullet$ 1) $(\Rightarrow)$. Let there exist at least one
isolated point $x_0 \in L \bigcap P$. Then for  $f$,  $f(x_0) \neq
0$ one has  $(f-f_i)(x_0)=f(x_0)$, thus the sequence $f_i$ does
not converge to $f$ with respect to the $\sup$ norm.

$(\Leftarrow)$. By assumption and construction of the
approximation process $\forall \varepsilon>0$ $\forall f \in
C_0(M)|_L$ $\exists f_i$ $\forall x \in \mathrm{supp} f$,
($((f_i-f)(x)<\varepsilon)$).

2) $(\Rightarrow)$. The algorithm finishes on the $i$-th step,
hence $d(x_j, x_k) \geq \frac{\varepsilon}{2^i}$. Since by
construction $(x_j)$ is compact it is finite.

$(\Leftarrow)$. Evident (cf. the third case of  \cite{Lob}).
$\triangleright$

Note here that the first case of \cite{Lob} can be described as
the set consisting only of the first component. The second case
then shows a possibility to find  set which does not contain the
first component of the decomposition. The third case is then the
case in which both components are present.

Let us consider a set $P' \subset P$ such that for some atlas
there exists $C>0$ $\|f^{n}\|<C$, $n \in \mathbb{N}$. Then we can
slightly deform a tail of the decomposition given above. It seems
natural that for some problems we shall need only those parts of
the decomposition whose Chesaro mean $n'/n$ is the maximal. In
certain cases we can prove the

\begin{note}
Consider two sets $S_1, S_2$ of points which define the second
part of the decomposition. Define for each $i \in (n_1)$ ($n_j$,
$j=1, 2$ is the set consisting of the intersection points of
$\mathbb{N}$ with $S_i$) $a_i$ which equals number of points from
$(n_2)$ that lie between the $i$-th and the $i+1$-st element from
$(n_1)$. If $a_i \leq K \in\mathbb{N}$ the the Hausdorff dimension
$d_1$ of the set $S_1 P$, defined as the limit $f^{(n_1)P}$ is not
greater than that defined as the limit $f^{(n_2)} P$.
\end{note}

$\bullet$ We must consider a covering $(A_i)$ of the set $S_1$. In
case there exists a set with the power non-equal to $0$ in the
mean such that for any $x_1 \in S_1$ there exists a number
$i(x_1)>0$ $i_2 \in [-i_1, i_1]$ ($f^{i_2}(x_1) \in S_2$ we can
construct a covering of the second set by the first one.
$\triangleright$

This statement holds true also in general case since for any
isolated point $x \in L \bigcap P$ $\dim\{x\}=0$.

\begin{example}{\em
Consider a foliation on $\mathbb{R}^2 \times [0, 1]/\simeq$, here
the relation $\simeq$ is given as follows: $(x, 0) \simeq (e^{2i}
\frac{x}{|x|^2}, 1)$. Assume that Ehresmann connection is given by
distribution tangent to $\mathbb{R}^2 \times\{pt\}$, $pt \in [0,
1]$. Then the maximal set for any leaf not passing through
$S_1=\{x \in \mathbb{R}^2| |x|=1\}$ is $S_1$ of dimension equal to
$1$, at the same time any other point of the intersection
naturally is isolated so of dimension $0$. Note that any leaf
passing through points inside $S_1$ has another limit set
--- point $\{0\}$ of dimension not less than that of any other point in the
intersection of this leaf and $P$.}
\end{example}

Note also that the Chesaro convergence of $n_1/n_2$ does not imply
the previous statement. It suffices to consider the first sequence
consisting of $(1/\sqrt{n} n)$. Assume that the members of the
second sequence equal to the members of the first one but their
number is such that $a_n=\sqrt{n}$.

Let us define a characteristic function for any point $x$ which is
one of the points generating the second part of the decomposition.
This function is built as follows: $\forall \varepsilon> 0$
consider the frequency of the intersection of the given leaf $L$
and $B_{\varepsilon}(x)$ $\nu_{\varepsilon}(x)=\lim\limits_{N \to
\infty}N_1/N$, here $N_1$ i number of the points on the leaf from
$[0, N]$ which lie in $B_{\varepsilon}(x)$. Consider
$\lim\limits_{\varepsilon \to 0} \nu_{\varepsilon}(x)$. Now the
most important function for any $\varepsilon>0$ is
$f_{\varepsilon}: \mathbb{N} \to \mathbb{N}$,
$f_{\varepsilon}=N_1(n)$. By construction $f_{\varepsilon} (n)
\leq k n$, $k \in [0, 1]$ and $f_{\varepsilon}> f_{\delta}$ for
$\delta < \varepsilon$. Set $f_0(x, k)=\sup\{n(k)| a^{n(k)} x_0
\to x\}$, ($x_0 \in L \bigcap P$ is some fixed point). This
mapping is similar to the return map of \cite{Milnor}, but is
greater than equal to the inverse of it.

Thus there either exists a function
$g(t)$, $g \neq \mathrm{const}$ $g(0)=0$ such that
$\lim\limits_{\varepsilon \to 0} \lim\limits_{t \to 0}
g(t)f_{\varepsilon}(1/t) = 1$ or $\lim\limits_{\varepsilon \to 0}
f_{\varepsilon}(1/t) = 0$. Let us show that in the last case $x$
is an isolated point. Let us consider $f(x, n)$ generated by the
sequence $x_n \in B_{1/n}(x) \bigcap L$. This function on infinity
is less or equal than any $f_{\varepsilon}(n)$ for any
$\varepsilon
>0$, at the same time by assumption $x$ is not isolated, thus
$\lim\limits_{n \to \infty}f(x, n)= \infty$. So the behavior of
the singular point is uniquely described by the behavior of the
function $g$ at zero and we can introduce an order on the set of
these points $x \geq y$ $\leftrightarrow$ $\lim\limits_{t \to
0}g(x, t)/g(y, t)< \infty$. Consider now $X_f=\{x \in P | f(x,
t)=f(t)\}$. Naturally the closure of this set contains only points
greater with respect to the relation given above. First mention
that on the one hand $\forall U_{\varepsilon}(x')$ $\exists
\varepsilon>0$, $x \in U_{\varepsilon'}(x')$,
($f_{\varepsilon'}(x')\geq f_{\varepsilon}(x)$). On the other hand
let $x_k \to y$ $a^{n_i+k} x_0 \to x_k$, then the sequence
$(x_j)$, $x_{2^k}, \ldots x_{2^{k+1}-1} \in (a^{n_i+k} x_0)_i$
$x_j \to y$ as $j \to \infty$ and the growth at $\infty$ of
$(n_i+k)_j$ coincides with that of $n_i$.

The same function can be defined for the set $S_0=\{p \in L\bigcap
P|$ there exists an open neighbourhood $U(p)$ of the point $p$
such that $U(p)\bigcap (L \bigcap P)=\{p\}\}$. Fist not that the set $S_0
\subset \mathbb{R}^n$ is at most countable one. Then  fix an order
on the set $\phi:\mathbb{N} \leftrightarrow L \bigcap P$. Note
that $\phi^{-1}(S_0)\subset \mathbb{N}$. Now for a point $s \in
S_0$, $s=\phi(0)$ we put $f_0(s)=0$. Finally for $\phi^{-1}(s) \in
\mathbb{Z}$ we get
$$
f(\phi^{-1}(s))=
\left\{%
\begin{array}{ll}
    f(\phi^{-1}(s-1))+1, & \hbox{if}\; s \in S_0;  \\
    f(\phi^{-1}(s-1)), & \hbox{otherwise.} \\
\end{array}%
\right.
$$

Now we point out that $X_M=\bigcap\limits_{f | X_f \neq \emptyset}
X_{g \leq f}$ is closed and non-empty in case $P$ or one of the
considered sets $X_{g \geq f}$ is compact. Thus we have proved the

\begin{statement}
If $P$ or one of the sets $X_{g \leq f}$ is compact then there
exists a maximal element on the set of points with the order
$\geq$.
\end{statement}

The statements of the following note are too evident to give the
proofs but rather important to ignore.

\begin{note}
{\em 1)} For any function $f$ the set $X_f$ is invariant with respect to
$H_P$.

{\em 2)} $X_{\min}=L \bigcap P$.

{\em 3)} $L \bigcap P \bigcap X_{\max} \neq \emptyset$ implies
$X_{\min}=X_{\max} = \overline{L \bigcap P}$.
\end{note}

\begin{corollary}
The function $f_0$ can either vanish or be equal to the trivial
function $f: \mathbb{N} \to \mathbb{N}$, $n \mapsto n$.
\end{corollary}

$\bullet$ To prove this consider a sequence of points $(x_n)$ of
$L: \bigcap P$ such that $x_n \to x \in L \bigcap P$. Now, $x
\not\in S_0$. Assume now that there also exists $x' \in L\bigcap
P$, $x' \in S_0$. Then there must exist $h \in H$ such that $x'=h
x$, but then since the action of the group $H$ is continuous $h
x_n \to x'$. This contradiction with the first statement of the
previous note proves the statement. $\triangleright$

\begin{example}{\em
Kronecker irrational flow on the torus $\mathbb{T}^2$. Consider
the parallel $P$ diffeomorphic to $\mathbb{S}^1$ of the torus
transversal to the leaves of the foliation. The set
$X_{\min}=X_{\max}=P$ since the only nonempty closed subset of
$\mathbb{S}^1$ invariant with respect to the rotation rationally
independent with $2 \pi$ is $\mathbb{S}^1$ itself. Note that in
case the set of points that return to the fixed one is the
sequence defined by Fibonacci numbers (appendix C of
\cite{Milnor}), then for any $x \in P$ $f_0(x, t) \geq
\log_{\frac{\sqrt{5}+1}{2}}\sqrt{5} t$. }\end{example}

Note that in case $\overline{L \bigcap P}\subset P$ does not
contain nontrivial sets invariant with respect to the action of $H_P$
then $X_{\max}=\overline{L \bigcap P}$. The following example
shows that the latter is possible also in other
case.

\begin{example}{\em Denjoit $C^1$-vector field on the torus $\mathbb{T}^2$
\cite{Tam}. Again consider the irrational flow on torus
$\mathbb{T}^2$ as a foliation. Consider next the countable set of
intervals $\{I_m=[0, l_m]; m \in \mathbb{Z}, l_m>0\}$ such that

1) $\sum\limits_{i \in \mathbb{Z}} l_i=l<+\infty$,

2) $\lim\limits_{m \to\infty} l_m/l_{m+1}=1$.

Substitute then the point $a^m x \in P$ by the interval $I_m$, $m
\in \mathbb{Z}$, where $P \cong \mathbb{S}^1$ is the parallel
transversal to the leaves of the foliation on the torus
$\mathbb{T}^2$, $a$ is the rotation of the transversal $P$ on the
angle rationally independent with $2 \pi$, $x \in P$ being a fixed
point.

Let us introduce now the set of mappings $f_m: I_m \to I_{m+1}$
with the following properties:

1) $\frac{d f_m}{d t}>0$;

2) there exists a number $\delta_m>0$ such that the derivative
$\frac{d f_m}{d t}$ equals $1$ on intervals $[0, \delta_m)$ and
$(l_m-\delta_m, l_m]$;

3) $\min(1, l_m/ l_{m+1})-(1-l_{m+1}/l_m)^2 \leq \frac{d f_m}{d t}
\leq \max(1, l_{m+1}/l_m)+(1-l_{m+1}/l_m)^2$.

Let us add this map to the shift $a: \mathbb{S}^1 \to
\mathbb{S}^1$. Thus we get a mapping $f: \mathbb{S}^1 \to
\mathbb{S}^1$.

Then the leaf passing through any point of $I_m$, $m \in
\mathbb{Z}$ intersects $P$ in closed nowhere dense set invariant
with respect to the map $f$. Thus the set of intersection points
of this leaf coincides with the $X_{\max}$ constructed above. Also
any leaf which does not pass through $I_m$, $m \in \mathbb{Z}$ is
dense on the set of similar leaves. Moreover the closure of the
set $L \bigcap P$ contains also points from $\partial(I_m)$. In
this case the maximal subset is a closure of $L \bigcap P \subset
P$ as in the preceding example. Note that this set contains the
invariant with respect to the action of $f$ closed subset
consisting of $\bigcup\limits_{m \in \mathbb{Z}}\partial I_m$. }
\end{example}

Since for any $g$ the set $X_g$ is invariant with respect to the
action of $H_P$ we can consider the following condition:

$(*)$ The characteristic function $f(x, t)$ does not depend on $x \in P$ in some neighbourhood
$U(x_0)$ of the fixed point $x_0  \in L \bigcap P$. The following
statement is then evident.

\begin{statement}
Condition $(*)$ implies $X_{\max} \subset \partial P \bigcap L$.
\end{statement}

\section{Family of Shr\"{o}dinger operators.}

Let us again consider a manifold $M$ with the foliation $F$
generated by the action of the abelian group $H$ and integrable
Ehresmann connection invariant with respect to the action of the
group $H$. Let us put into consideration the operator
$\mathcal{H}=-\nabla+V$ on $L^2(H \times P)$, defined for any leaf
of the foliation $F$, here $V \in C_0(M)$ and
$\nabla=\sum\limits_{i=1}^{\dim H} \frac{\partial^2}{\partial
x_i^2} $, $x_i$ $i=\overline{1, \dim H}$ being coordinates on $H$. Thus $\mathcal{H}$
defines a family of Shr\"{o}dinger operators $(\mathcal{H}_p)_{p
\in H}$, $\mathcal{H}_p: L^2(\mathbb{R}^n) \to L^2(\mathbb{R}^n)$
for any leaf $L$, here $n=\dim H$. So each
$\mathcal{H}_p=-\nabla+V_p$ for $V_p \in C_0(M)|_{L_p}$. At the
rest of the paper we assume that $V_p$ depends upon $p \in P$
analytically. The good illustration of the method gives the paper \cite{LP}.

\begin{statement}
Let $H=\mathbb{S}^1$. Then the spectrum of the Shr\"odinger
operator $\mathcal{H}=-\frac{d^2}{d x^2}+V$ depends on the transversal
coordinate $p \in P$ continuously.
\end{statement}

$\bullet$ Let $\mu$ be an $H$-invariant measure on the transversal
$P$. Consider the decomposition of the set of the target spaces
$L_2$
$$
L^2(H
\times P)=\int\limits_{P}L^2(H) d \mu =\int\limits_{P}
\int\limits_{[0, 2 \pi)^n}H'\frac{d^n \theta}{(2 \pi)^n} d
\mu,
$$
here $H'\simeq l_2$ and the isomorphism  $L^2(\mathbb{R}^n = H,
d^n x) \leftrightarrow \int\limits_{[0, 2 \pi)^n}H'\frac{d^n
\theta}{(2 \pi)^n}$ is the Floquet decomposition $U':
L^2(\mathbb{R}^n = H, d^n x) \to \int\limits_{[0, 2
\pi)^n}H'\frac{d^n \theta}{(2 \pi)^n}$ given as follows:
$$
(U'f)_{\theta}(x)=\sum\limits_{m \in
\mathbb{Z}^n} e^{- i \theta m} f(x +\sum m_i a_i),
$$
$a_i$, $i=\overline{1, n}$ being the basis of the group $H$ (not
such of the group $H_P$). Then theorems XIII.64 and XII.11 from
\cite{RS4} and Foubini theorem imply the statement.
$\triangleright$

Assume now that $H$ is not a compact group and for any $x \in P$
the set $H_P\{x\}$ is dense in $P$. Assume also that there exists
a $H_P$-invariant metric on $P$. Thus we deal with specific
ergodic transformations on $P$. Then it can be shown that
$C_0(M)|_{L}$ consists of limit almost-periodic functions with
periods rational dependent with $a \in H_P$.

\begin{lemma}
The spectrum of the Shr\"odinger operator $\mathcal{H}=-\frac{d}{d
x}+V$ with limit almost periodic potential $V$ is a pointwise
limit of a set of spectra of operators $\mathcal{H}_n=-\frac{d}{d
x}+V_n$, here $V_n$ are limit almost-periodic functions such that
$V_n \to_{\|\cdot\|_{\sup}} V$, $(n \to \infty)$.
\end{lemma}

$\bullet$ Consider $\varepsilon>0$ and $\|\mathcal{H}_n -
\mathcal{H}\|\leq \varepsilon$. Let us also take into
consideration spaces $L_{n, m}= \int\limits{(-1/2, 1/2]} l' d q$,
where $l'=l_2(\{1, \ldots, m\})$ and the operator $U$ is again a
Floquet decomposition. Then one have for the set of operators
$\mathcal{H}^{'}_n =\mathcal{H}_n|_{L_{n, m}}$. At the same time
$\mathcal{H}_n$ is a bounded operator for any $n \in \mathbb{N}$.
Then Corollary 2 of \cite{DSch} implies the statement for
$\mathcal{H}^{'}_n$, thus any point of the spectrum of
$\mathcal{H}$ is close to some such point of $\mathcal{H}_n$ for
sufficiently large number $n \in \mathbb{N}$. $\triangleright$

\begin{corollary}
If there exists an $H_P$-invariant riemannian metric on $P$ then
spectrum of the family of Schr\"odinger operators continuously depend on $p
\in P$.
\end{corollary}

\begin{note}
Note that similar statements can be translated to the case of the
complex-valued potential $V : M \to \mathbb{C}$. This can be done
using results of \cite{RB} \cite{Shin}.
\end{note}

\begin{statement}
Let transversal $P \simeq \mathbb{R}^n$ and $\forall x \neq 0$
$H_x = n \mathbb{Z}$ ($n \in \mathbb{N}$), $H_{0}=\mathbb{Z}$.
Then the following statements hold true:

{\em 1)} Spectrum of the family of Hill operators with periodic
real-valued potentials is homeomorphic to $L \times
\mathbb{S}^{n-1}$, here $L$ is the union of graphs of functions
$f_{i}: \mathbb{R} \to \mathbb{R}$ $i \in \mathbb{N}$, $f_{n
i+k}(x)= 1/k \arctan(x)$, $f_{n i}=0$, $k = \overline{1, \ldots,
n-1}$, $i \in \mathbb{N}$.

{\em 2)} Assume that for any $x \in P$ and $E \in \Lambda_{x}$
($\Lambda_x$ is the spectrum of the operator $H_p$ at the point $x
\in P$) $\Delta^{'}_x(E) \neq 0$. If $V$ is a complex-valued
function then arcs of the spectrum with the greatest period can
interchange rotating over $0 \in P$ only in case operator $H_0$ at
the point $0$ possess closed analytic arc in the spectrum.
\end{statement}

$\bullet$ $V$ is the real-valued function. There exists only one
non-trivial case --- $P \simeq \mathbb{R}^2$, since the only
alternative to the given case is the spiral line over $\{0\} \times
\mathbb{R} \in P \times \mathbb{R}$. Then the spectrum of each
operator $H_{\theta}$ is not purely discreet one which contradicts
paragraph $(a)$ of the Theorem XIII.89 from \cite{RS4}.

Assume that the condition of the paragraph 2) holds true then
passing to the limit $p \to 0$ we obtain the arc in the spectrum
of $H_0$ analytic in each internal point \cite{Shin}. Then by
periodicity of the potential $V$ the arc of the operator $H_0$,
which by the last note splits into $n$ arcs of the operators
$H_x$, $x \in P \setminus\{0\}$ is the connected set i.e. $n-1$
arcs of the spectrum of the operator $H_x$ for $x \neq 0$ glue
with each other. Thus invariance of these arcs under the rotation
over $0$ --- which at the same time gives rise to the
transformation of the spaces $L^2(\mathbb{R})$ --- the other arcs
also glue. $\triangleright$

Using the density functions of the second part of the paper one can
find out which eigenvalues of the family of operators are more
important.


\begin{thebibliography}{99}
\bibitem{BH1} R.A. Blumenthal and J.J. Hebda,
An analogue of the holonomy bundle for a foliated manifold.
C. R. Acad. Sci., Paris, S\'er. I 303  (1986), 931 -- 934.

\bibitem{BH} R.A. Blumenthal and J.J. Hebda, Complementary
distributions which preserve the
leaf geometry
and applications to totally geodesic foliations. Quart. J. Math. Oxford 35 (1984), 383 -- 392.

\bibitem{BH2} R.A. Blumenthal and J.J. Hebda,
Ehresmann connections for foliations.
Indiana Univ. Math. J. 33  (1984), 597 -- 611.

\bibitem{Cadet} F. Cadet, Deformation quantization using groupoids.
Case of toric manifolds, arXiv:math.OA/0305261 v2 20 May 2003.

\bibitem{DSch} N. Dunford and  J.T. Schwartz, Linear Operators. I.
General theory. (Pure and Applied Mathematics. Vol. 6). New York and London:
Interscience Publishers. 1958.

\bibitem{Lob}  P. N. Ivanshin,
Structure of function algebras on foliated manifolds.
Lobachevskii J. Math. 14 electronic only (2004), 39 -- 54.

\bibitem{SW}  P. N. Ivanshin,
Algebras of functions on groupoid of some special foliations.
Southwest J. Pure Appl. Math. 2003 No.1, electronic only (2003), 96 -- 108.

\bibitem{KoNo} Sh. Kobayashi and K. Nomizu,
Foundations of differential geometry. I. New York-London:
Interscience Publishers, a division of John Wiley \& Sons. 1963.

\bibitem{LP} F. \'Lledo and O. Post, Generating spectral
gaps by geometry, arXiv:math-ph/0406032 v1 15 Jun 2004.

\bibitem{Milnor} J. Milnor, Dynamics in one complex
variable. Introductory lectures. Wiesbaden: Vieweg. 1999.

\bibitem{Mo} P. Molino,
Riemannian foliations.
Progress in Mathematics, Vol. 73. Boston-Basel: Birkh\"auser. 1988.

\bibitem{Nov} S.P. Novikov,
Topology of foliations. (Russian, English) Trans. Mosc. Math. Soc.
14 (1965), 268 -- 304;  translation from Tr. Mosk. Mat. Obshch. 14
(1965), 248 -- 278.

\bibitem{RS4} M. Reed and B. Simon,
Methods of modern mathematical physics. IV: Analysis of operators.
New York - San Francisco - London: Academic Press. 1978.

\bibitem{RB} F.S. Rofe-Beketov,
The spectrum of non-selfadjoint differential operators with
periodic coefficients (Russian, English) Sov. Math., Dokl. 4
(1963), 1312 -- 1315; translation from Dokl. Akad. Nauk SSSR 152
(1963), 1312-1315 .

\bibitem{Shin} K.C. Shin,
On the shape  of spectra for non-self-adjoint periodic
Schr\"odinger operators. arXiv:math-ph/0404015 v1 6 Apr 2004.

\bibitem{Tam} I. Tamura, Topology
of foliations: an introduction. Translations of Mathematical
Monographs. 97. Providence, RI: American Mathematical Society
(AMS). 1992.

\bibitem{ZhCh} N. I. Zhukova and  G.V. Chubarov,
Aspects of the qualitative theory of suspended foliations.
J. Difference Equ. Appl. 9, No. 3-4 (2003), 393 -- 405.

\end{thebibliography}
\end{document}